\documentclass[11pt]{article}
\usepackage{amsmath}
\usepackage{amssymb}
\usepackage{latexsym}
\usepackage{dsfont}
\def\disp{\displaystyle}

\def\F{Fr\'{e}chet}

\def\tto{\;{\lower 1pt \hbox{$\rightarrow$}}\kern -10pt
\hbox{\raise 2pt \hbox{$\rightarrow$}}\;}
\def\lto{\longrightarrow}
\def\Hat{\widehat}
\def\Tilde{\widetilde}

\def\ra{\rangle}
\def\la{\langle}
\def\ve{\varepsilon}
\def\B{I\!\!B}
\def\h{\hfill\Box}
\def\R{\mathbb{R}}
\def\N{\mathbb{N}}

\def\ox{\bar{x}}

\def\co{\mbox{\rm co}\,}

\def\epi{\mbox{\rm epi}\,}

\def\cl{\mbox{\rm cl}\,}
\def\dist{\mbox{\rm dist}\,}
\def\cone{\mbox{\rm cone}\,}

\def\inter{\mbox{\rm int}\,}

\def\mT{\mathcal{T}}

\def\h{\hfill\square}
\def\dn{\downarrow}
\def\O{\Omega}

\def\ph{\varphi}
\def\emp{\emptyset}
\def\st{\stackrel}

\def\lm{\lambda}
\def\Lm{\Lambda}
\def\gg{\gamma}

\def\dd{\delta}
\def\al{\alpha}

\def \N{I\!\!N}

\begin{document}
\begin{center}
{\bf TANGENTIAL EXTREMAL PRINCIPLES FOR FINITE AND INFINITE SYSTEMS OF SETS,
I: BASIC THEORY}\footnote{This research was partially supported by the US National Science
Foundation under grants DMS-0603846 and DMS-1007132 and by the Australian Research Council
under grant DP-12092508.}\\[2ex]
BORIS S. MORDUKHOVICH\footnote{Department of Mathematics, Wayne
State University, Detroit, MI 48202, USA. Email: boris@math.wayne.edu.} and HUNG M. PHAN\footnote{Department of Mathematics, Wayne
State University, Detroit, MI 48202, USA. Email:
pmhung@wayne.edu.}
\end{center}
\newtheorem{Theorem}{Theorem}[section]
\newtheorem{Proposition}[Theorem]{Proposition}
\newtheorem{Lemma}[Theorem]{Lemma}
\newtheorem{Corollary}[Theorem]{Corollary}
\newtheorem{Definition}[Theorem]{Definition}
\newtheorem{Remark}[Theorem]{Remark}
\newtheorem{Example}[Theorem]{Example}
\renewcommand{\theequation}{\thesection.\arabic{equation}}
\normalsize
\setcounter{equation}{0}
\newcounter{count}
\newenvironment{lista}{\begin{list}{(\rm\alph{count})}{\usecounter{count}
\setlength{\itemindent}{-0.3cm}
}}{\end{list}}
\newenvironment{listi}{\begin{list}{\rm(\textit{\roman{count}})}{\usecounter{count}
\setlength{\itemindent}{-0.3cm}}}{\end{list}}

\small{{\bf Abstract.} In this paper we develop new extremal principles in
variational analysis that deal with finite and infinite systems of convex and
nonconvex sets. The results obtained, unified under the name of tangential
extremal principles, combine primal and dual approaches to the study of variational
systems being in fact first extremal principles applied to infinite systems of sets.
The first part of the paper concerns the basic theory of tangential extremal principles
while the second part presents applications to problems of semi-infinite programming
and multiobjective optimization.
\vspace*{0.1cm}

\vspace*{0.1cm}

}
\normalsize

\section{Introduction}
\setcounter{equation}{0}

It has been well recognized that the {\em convex separation
principle} plays a crucial role in many aspects of nonlinear
analysis, optimization, and their applications. In particular, a
conventional way to derive necessary optimality conditions in constrained
optimization problems is to construct first an appropriate {\em tangential convex
approximations} of the problem data around an optimal solution in
{\em primal} spaces and then to apply a convex separation theorem to get
supporting elements in dual spaces (Lagrange multipliers, adjoint
arcs, shadow prices, etc.). For problems of nonsmooth optimization, this approach
inevitably leads to the usage of {\em convex} sets of normals and
subgradients whose calculi are also based on convex separation
theorems and/or their equivalents.

Despite the well-developed technique of convex analysis, the convex separation approach
has a number of serious limitations, especially concerning applications to problems of nonsmooth
optimization and related topics; see, e.g., commentaries and discussions on pp.\ 132--140 of
\cite{m-book1} and also on pp.\ 131--133 of \cite{m-book2}. To overcome some of these limitations,
a {\em dual-space} approach revolving around {\em extremal principles} has been developed and
largely applied in the frameworks of variational analysis, generalized differentiation, and
optimization-related areas; see the two-volume monograph \cite{m-book1,m-book2} with their
references. The extremal principles developed therein can be viewed as variational counterparts of
convex separation theorems in nonconvex settings while providing normal cone descriptions of extremal
points of finitely many closed sets in terms of the corresponding generalized {\em Euler equation}.

Note that the known extremal principles do not involve any tangential approximations of sets in
primal spaces and do not employ convex separation. This dual-space approach exhibits a number of
significant advantages in comparison with convex separation techniques and opens new perspectives in
variational analysis, generalized differentiation, and their numerous applications. On the other hand,
we are not familiar with any versions of extremal principles in the scope of \cite{m-book1,m-book2}
for infinite systems of sets; it is not even clear how to formulate them appropriately in the lines of the
developed methodology. Among primary motivations for considering infinite systems of sets we mention problems
of semi-infinite programming, especially those concerning the most difficult case of countably many constraints
vs.\ conventional ones with compact indexes; cf.\ \cite{gl}.\vspace*{0.05in}

The main purpose of this paper is to propose and justify extremal principles of a new type, which can be
applied to infinite set systems while also provide independent results for finitely many nonconvex sets.
To achieve this goal, we develop a novel approach that incorporates and unifies some ideas from both
tangential approximations of sets in primal spaces and nonconvex normal cone approximations in dual spaces.
The essence of this approach is as follows. Employing a variational technique, we first derive a new
{\em conic extremal principle}, which concerns countable systems of general {\em nonconvex} cones in
finite dimensions and describes their extremality at the origin via an appropriate countable version of the
generalized Euler equation formulated in terms of the nonconvex limiting normal cone by Mordukhovich
\cite{Mor76}. Then we introduce a notion of {\em tangential extremal points} for infinite (in particular, finite)
systems of closed sets involving their tangential approximations. The corresponding {\em tangential extremal
principles} are induced in this way by applying the conic extremal principle to the collection of selected tangential
approximations. The major attention is paid in this paper to the case of tangential approximations generated by
the (nonconvex) Bouligand-Severi {\em contingent cone}, which exhibits remarkable properties that are most
appropriate for implementing the proposed scheme and subsequent applications. The contingent cone is replaced by its
{\em weak} counterpart when the space in question is infinite-dimensional. Selected applications of the developed theory
to problems of semi-infinite programming and multiobjective optimization are given in the second part of this study
\cite{m-hung}\vspace*{0.05in}

For the reader's convenience we briefly overview in Section~2 some basic constructions of tangent and normal cones in
variational analysis widely used in what follows. Section~3 contains definitions of tangential extremal points
of finite and infinite set systems as well as descriptions of the extremality conditions for them, which are at
the heart of the tangential extremal principles established below. In this section we also compare the new notions of
tangential extremality with the conventional notion of extremality previously known for finite systems of sets.

Section~4 is devoted to deriving the conic extremal principle for countable systems of arbitrary closed cones in
finite-dimensional spaces. In Section~5 we apply this basic result to establishing several useful representations of
Fr\'echet normals to countable intersections of cones at the origin.

Section~6 concerns the study of the weak contingent cone in
infinite-dimensional spaces, which reduces to the classical
Bouligand-Severi contingent cone in finite dimensions. We show
that the weak contingent cone provides a remarkable tangential
approximation for an arbitrary closed subset enjoying, in
particular, the new tangential normal enclosedness and approximate
normality properties in any reflexive Banach spaces. These
properties are employed In Section~7 to derive contingent and weak
contingent extremal principles for countable and finite systems of
closed sets in finite and infinite dimensions. We also establish
appropriate versions of the aforementioned results in a broader
class of Asplund spaces. \vspace*{0.05in}

Throughout the paper we use standard notation of variational
analysis; see, e.g., \cite{m-book1,Rockafellar-Wets-VA}. Unless
otherwise stated, the space $X$ in question is Banach with the
norm $\|\cdot\|$ and the canonical pairing $\la\cdot,\cdot\ra$
between $X$ and its topological dual $X^*$ with $\B\subset X$ and
$\B^*\subset X^*$ standing for the corresponding closed unit
balls. The symbols $\st{w}{\to}$ and $\st{w^*}{\to}$ indicate the
weak convergence in $X$ and the weak$^*$ convergence in $X^*$,
respectively. Given $\emp\ne\O\subset X$, denote by
\begin{equation*}
\cone\O:=\bigcup_{\lm\ge 0}\lm\O=\bigcup_{\lm\ge 0}\Big\{\lm v\Big|\;v\in\O\Big\}
\end{equation*}
the {\em conic hull} of $\O$ and by
\begin{equation*}
\co\O:=\Big\{\sum_{i\in I}\lm_i u_i\Big|\;I\;\mbox{ finite },\;\lm_i\ge 0,\;\sum_{i\in I}\lm_i=1,\;u_i\in\O\Big\}
\end{equation*}
the {\em convex hull} of this set. The notation $x\st{\O}{\to}\ox$
means that $x\to\ox$ with $x\in\O$.  Finally, $\N:=\{1,2,\ldots\}$
signifies the collection of all natural numbers.

\section{Tangents and Normal to Nonconvex Sets}
\setcounter{equation}{0}

In this section we recall some basic notions of tangent and normal
cones to nonempty sets closed around the reference points; see the
books \cite{Borwein-Zhu-TVA,m-book1,Rockafellar-Wets-VA,s} for
more details and related material.

Given $\O\subset X$ and $\ox\in\O$, the closed (while often nonconvex) cone
\begin{equation}\label{bs}
T(\ox;\O):=\big\{v\in X\big|\;\exists\,\mbox{ sequences }\;t_k\dn
0,\;v_k\to v\;\mbox{ with }\;\ox+t_k v_k\in\O,\;\forall
k\in\N\big\}
\end{equation}
is the {\em Bouligand-Severi tangent/contingent cone} to $\O$ at $\ox$. We also use its weak counterpart
\begin{equation}\label{wbs}
T_w(\ox;\O):=\big\{v\in X\big|\;\exists\,\mbox{ sequences
}\;t_k\dn 0,\;v_k\st{w}{\to}v\;\mbox{ with }\;\ox+t_k v_k\in\O,\;
\forall k\in\N\big\}
\end{equation}
known as the {\em weak contingent cone} to $\O$ at this point. For any $\ve\ge 0$, the collection
\begin{equation}\label{eq:e-nor}
\Hat N_\ve(\ox;\O):=\left\{x^*\in
X^*\Big|\limsup_{x\st{\O}{\to}\ox}\frac{\la
x^*,x-\ox\ra}{\|x-\ox\|}\le\ve\right\}
\end{equation}
is called the set of {\em $\ve$-normals} to $\O$ at $\ox$. In the
case of $\ve=0$ the set $\Hat N(\ox;\O):=\Hat N_0(\ox;\O)$ is a
cone known as the {\em Fr\'echet/regular normal cone} (or the
prenormal cone) to $\O$ at this point. Note that the Fr\'echet
normal cone is always convex while it may be trivial (i.e.,
reduced to $\{0\}$) at boundary points of simple nonconvex sets in
finite dimensions as for $\O=\{(x_1,x_2)\in\R^2|\;x_2
\ge-|x_1|\big\}$ at $\ox=(0,0)$. If the space $X$ is reflexive,
then
\begin{equation}\label{du}
\Hat N(\ox;\O)=T^*_w(\ox;\O):=\big\{x^*\in X^*\big|\;\la x^*,v\ra\le 0,\;\forall v\in T_w(\ox;\O)\big\}.
\end{equation}
The collection of sequential limiting normals
\begin{equation}\label{nc}
\begin{array}{ll}
N(\ox;\O):=\Big\{x^*\in X^*\Big|&\exists\;\mbox{ sequences }\;\ve_k\dn 0,\;x_k\st{\O}{\to}\ox,\;x^*_k\st{w^*}{\to}x^*\;
\mbox{ as }\;k\to\infty\\
&\mbox{such that }\;x^*_k\in\Hat N_{\ve_k}(x_k;\O),\,\forall k\in\N\Big\}
\end{array}
\end{equation}
is known as the {\em Mordukhovich/basic/limiting normal cone} to
$\O$ at $\ox$. If the space $X$ is Asplund, i.e., each of its
separable subspaces has a separable dual (this is automatic, in
particular, for any reflexive Banach space), then we can
equivalently put $\ve_k=0$ in \eqref{nc}; see \cite{m-book1} for
more details. Observe also that for $X=\R^n$ the normal cone
\eqref{nc} can be equivalently described in the form
\begin{equation*}
\begin{array}{ll}
N(\ox;\O)=\Big\{x^*\in\R^n\Big|&\exists\,\mbox{ sequences }\;x_k\to\ox,\;w_k\in\Pi(x_k;\O),\;\al_k\ge 0\\
&\mbox{such that }\;\al_k(x_k-w_k)\to x^*\;\mbox{ as }\;k\to\infty\Big\}
\end{array}
\end{equation*}
via the {\em Euclidean projector} $\Pi(x;\O):=\{w\in\O|\;\|x-w\|=\mbox{dist}(x;\O)\}$ of $x\in\R^n$ onto $\O$.

It is worth mentioning that the limiting normal cone \eqref{nc} is
often nonconvex as, e.g., for the set $\O\subset\R^2$ considered
above, where $N(0;\O)=\{(u_1,u_2)\in\R^2|\;u_2=-|u_1|\}$. It does
not happen when $\O$ is {\em normally regular} at $\ox$ in the
sense that $N(\ox;\O)=\Hat N(\ox;\O)$. The latter class includes
convex sets when both cones \eqref{eq:e-nor} as $\ve=0$ and
\eqref{nc} reduce to the classical normal cone of convex analysis
and also some other collections of ``nice" sets of a certain
locally convex type. At the same time it excludes a number of
important settings that frequently appear in applications; see,
e.g., the books \cite{m-book1,m-book2,Rockafellar-Wets-VA} for
precise results and discussions. Being nonconvex, the normal cone
$N(\ox;\O)$ in \eqref{nc} cannot be tangentially generated by
duality of type \eqref{du}, since the duality/polarity operation
automatically implies convexity. Nevertheless, in contrast to
Fr\'echet normals, this limiting normal cone enjoys {\em full
calculus} in general Asplund spaces, which is mainly based on
extremal principles of variational analysis and related
variational techniques; see \cite{m-book1} for a comprehensive
calculus account and further references. \vspace*{0.05in}

The next simple observation is useful in what follows.

\begin{Proposition}\label{Prop:N(Lm)} {\bf (generalized normals to cones).} Let $\Lm\subset X$ be a cone, and
let $w\in\Lm$. Then we have the inclusion
\begin{equation*}
\Hat N(w;\Lm)\subset N(0;\Lm).
\end{equation*}
\end{Proposition}
{\bf Proof.} Pick any $x^*\in\Hat N(w;\Lm)$ and get by definition \eqref{eq:e-nor} of  the \F\
normal cone that
\begin{equation*}
\limsup_{x\st{\Lm}{\to}w}\frac{\la x^*,x-w\ra}{\|x-w\|}\leq 0.
\end{equation*}
Fix $x\in\Lm$, $t>0$ and let $u:=x/t$. Then $(x/t)\in\Lm$, $tw\in\Lm$, and
\begin{equation*}
\limsup_{x\st{\Lm}{\to}tw}\frac{\la x^*,x-tw\ra}{\|x-tw\|}=
\limsup_{x\st{\Lm}{\to}w}\frac{t\la
x^*,(x/t)-w\ra}{t\|(x/t)-w\|}=
\limsup_{u\st{\Lm}{\to}w}\frac{\la x^*,u-w\ra}{\|u-w\|}\le0,
\end{equation*}
which gives $x^*\in\Hat N(tw;\Lm)$ by \eqref{eq:e-nor}. Letting finally  $t\to 0$, we get $x^*\in
N(0;\Lm)$ and thus complete the proof of the proposition. $\h$

\section{Tangential Extremal Systems and Extremality Conditions}
\setcounter{equation}{0}

In this section we introduce the notions of conic and tangential
extremal systems for finite and countable collections of sets and
discuss extremality conditions, which are at the heart of the
conic and tangential extremal principles justified in the
subsequent sections. These new extremality concepts are compared
with conventional notions of local extremality for set systems.
\vspace*{0.05in}

We start with the new definitions of extremal points and extremal
systems of a countable or finite number of cones and general sets
in normed spaces.

\begin{Definition}\label{ext-def} {\bf (conic and tangential extremal systems).} Let $X$ be an arbitrary normed space.
Then we say that:

{\bf (a)} A countable system of cones $\{\Lm_i\}_{i\in\N}\subset X$ with $0\in\cap_{i=1}^\infty\Lm_i$ is
{\sc extremal at the origin}, or simply is an {\sc extremal system of cones}, if there is a bounded sequence
$\{a_i\}_{i\in\N}\subset X$ with
\begin{equation}\label{cone-emp}
\bigcap_{i=1}^\infty\big(\Lm_i-a_i\big)=\emp.
\end{equation}

{\bf (b)} Let $\{\O_i\}_{i\in\N}\subset X$ be an countable system
of sets with $\ox\in\cap_{i=1}^\infty\O_i$, and let
$\Lm:=\{\Lm_i(\ox)\}_{i\in\N}$ with
$0\in\cap_{i=0}^\infty\Lm_i(\ox)\subset X$ be an approximating
system of cones. Then $\ox$ is a $\Lm$-{\sc tangential local
extremal point} of $\{\O_i\}_{i\in\N}$ if the system of cones
$\{\Lm_i(\ox)\}_ {i\in\N}$ is extremal at the origin. In this case
the collection $\{\O_i,\ox\}_{i\in\N}$ is called a $\Lm$-{\sc
tangential extremal system}.

{\bf (c)} Suppose that $\Lm_i(\ox)=T(\ox;\O_i)$ are the contingent
cones to $\O_i$ at $\ox$ in \rm(b). Then $\{\O_i,\ox\}_{i\in\N}$
is called a {\sc contingent extremal system} with the {\sc
contingent local extremal point} $\ox$. We use the terminology of
{\sc weak contingent extremal system} and {\sc weak contingent
local extremal point} if $\Lm_i(\ox)=T_w(\ox;\O_i)$ are the weak
contingent cones to $\O_i$ at $\ox$.
\end{Definition}

Note that all the notions in Definition~\ref{ext-def} obviously
apply to the case of systems containing {\em finitely} many sets;
indeed, in such a case the other sets reduce to the whole space
$X$. Observe also that both parts in part (c) of this definition
are equivalent in finite dimensions. Furthermore, they both reduce
to (a) in the general case if all the sets $\O_i$ are cones and
$\ox=0$. \vspace*{0.05in}

Let us now compare the new notions of Definition~\ref{ext-def}
with the conventional notion of locally extremal points for
finitely many sets first formulated in \cite{KruMor80}. Recall
\cite[Definition~2.1]{m-book1} that a point
$\ox\in\cap_{i=1}^m\O_i$ is {\em locally extremal} for the system
$\{\O_1,\ldots,\O_m\}$ if there are sequences $\{a_{ik}\}\subset
X$ with $a_{ik}\to 0$ as $k\to\infty$ for $i=1,\ldots,m$ and a
neighborhood $U$ of $\ox$ such that
\begin{equation}\label{eq:ESfin}
\bigcap_{i=1}^m\big(\O_i-a_{ik}\big)\cap U=\emp\;\mbox{ for all large }\;k\in\N.
\end{equation}

We first observe that for finite systems of cones the local
extremality of the origin in the sense of \eqref{eq:ESfin} is
equivalent to the validity of condition \eqref{cone-emp} of
Definition~\ref{ext-def}.

\begin{Proposition}\label{cone-eq} {\bf (equivalent description of cone extremality).} The finite system of cones
$\{\Lm_1,\ldots,\Lm_m\}$ is extremal at the origin in the sense of Definition~{\rm\ref{ext-def}(a)} if and only if
$\ox=0$ is a local extremal point of $\{\Lm_1,\ldots,\Lm_m\}$ in the sense of \eqref{eq:ESfin}.
\end{Proposition}
{\bf Proof.} The ``only if" part is obvious. To justify the ``if"
part, assume that there are elements $a_1,\ldots,a_m\in X$ such
that
\begin{equation}\label{cone-eq1}
\bigcap_{i=1}^{m}\big(\Lm_i-a_i\big)=\emp.
\end{equation}
Now for any $\eta>0$ we have by \eqref{cone-eq1} and the conic structure of $\Lm_i$ that
\begin{equation*}
\emp=\bigcap_{i=1}^{m}\eta\big(\Lm_i-a_i\big)=\bigcap_{i=1}^{m}\big(\eta\Lm_i-\eta a_i\big)=\bigcap_{i=1}^{m}
\big(\Lm_i-\eta a_i\big).
\end{equation*}
Letting $\eta\dn0$ implies that the extremality condition
\eqref{eq:ESfin} holds, i.e., the origin is a local extremal point
of the cone system $\{\Lm_1,\ldots,\Lm_m\}$. $\h$\vspace*{0.05in}

Next we show that the local extremality \eqref{eq:ESfin} and the
contingent extremality from Definition~\ref{ext-def}(c) are
independent notions even in the case of two sets in $\R^2$.

\begin{Example}\label{exam-ext} {\bf (contingent extremality versus local extremality).}

{\rm {\bf(i)} Consider two closed subsets in $\R^2$ defined by
\begin{equation*}
\O_1:=\epi\ph\;\mbox{ with }\;\ph(x):=x\sin(1/x)\;\mbox{ as
}\;x\ne 0,\;\ph(0)=0\;\mbox{ and
}\;\;\O_2:=(\R\times\R_-)\setminus{\rm int}\,\O_1.
\end{equation*}
Take the point $\ox=(0,0)\in\O_1\cap\O_2$ and observe that the
contingent cones to $\O_1$ and $\O_2$ at $\ox$ are computed,
respectively, by
\begin{equation*}
T(\ox;\O_1)=\epi(-|\cdot|)\;\mbox{ and }\;T(\ox;\O_2)=\R\times\R_-.
\end{equation*}
It is easy to see that $\ox$ is a local extremal point of
$\{\O_1,\O_2\}$ but not a contingent local extremal point of this
set system.

{\bf(ii)} Define two closed subsets of $\R^2$ by
\begin{equation*}
\O_1:=\big\{(x_1,x_2)\in\R^2\big|\;x_2\ge-x^2_1\big\}\;\mbox{ and }\;\O_2:=\R\times\R_-.
\end{equation*}
The contingent cones to $\O_1$ and $\O_2$ at $\ox=(0,0)$ are computed by
\begin{equation*}
T(\ox;\O_1)=\R\times\R_+\;\mbox{ and }\;T(\ox;\O_2)=\R\times\R_-.
\end{equation*}
We can see that $\{\O_1,\O_2,\ox\}$ is a contingent extremal system but not an extremal system of sets.}
\end{Example}

Our further intention is to derive verifiable {\em extremality
conditions} for tangentially extremal points of set systems in
certain {\em countable} forms of the {\em generalized  Euler
equation} expressed via the limiting normal cone \eqref{nc} at the
points in question. Let us first formulate and discuss the desired
conditions, which reflect the essence of the tangential extremal
principles of this paper.

\begin{Definition}\label{Def:EPinf} {\bf(extremality conditions for countable systems).} We say that:

{\bf (a)} The system of cones $\{\Lm_i\}_{i\in\N}$ in $X$
satisfies the {\sc conic extremality conditions} at the origin if
there are normals $x^*_i\in N(0;\Lm_i)$ for $i=1,2,\ldots$ such
that
\begin{equation}\label{eq:EPinf}
\sum_{i=1}^{\infty}\frac{1}{2^i}x^*_i=0\quad\mbox{and}\quad\sum_{i=1}^{\infty}\frac{1}{2^i}\|x^*_i\|^2=1.
\end{equation}

{\bf (b)} Let $\{\O_i\}_{i\in\in\N }$ with
$\ox\in\cap_{i=1}^\infty\O_i$ and $\Lm:=\{\Lm_i\}_{i\in\N}$ with
$0\in\cap_{i=1}^\infty\Lm_i$ be, respectively, systems of
arbitrary sets and approximating cones in $X$. Then the system
$\{\O_i\}_{i\in\N}$ satisfies the $\Lm$-{\sc tangential
extremality conditions} at $\ox$ if the systems of cones
$\{\Lm_i\}_{i\in\N}$ satisfies the conic extremality conditions at
the origin. We specify the {\sc contingent extremality conditions}
and the {\sc weak contingent extremality conditions} for
$\{\O_i\}_{i\in\N}$ at $\ox$ if $\Lm=\{T(\ox;\O_i)\}_{i\in\N}$ and
$\Lm=\{T_w(\ox;\O_i)\}_{i\in\N}$, respectively.

{\bf (c)} The system of sets $\{\O_i\}_{i\in\N}$ in $X$ satisfies
the {\sc limiting extremality conditions} at $\ox\in
\cap_{i=1}^\infty\O_i$ if there are limiting normals $x^*_i\in
N(\ox;\O_i)$, $i=1,2,\ldots$, satisfying \eqref{eq:EPinf}.
\end{Definition}

Let us briefly discuss the introduced extremality conditions.

\begin{Remark} \label{ext-dis} {\bf (discussions on extremality conditions).}

{\rm{\bf (i)} All the conditions of Definition~\ref{Def:EPinf} can
be obviously specified to the case of {\em finite systems} of sets
by considering all the other sets as the whole space therein. Then
the series in \eqref{eq:EPinf} become finite sums and the
coefficients $2^{-i}$ can be dropped by rescaling.

{\bf (ii)} It easily follows from the constructions involved that
the contingent, weak contingent, and limiting extremality
conditions are are {\em equivalent} to each other if all the sets
$\O_i$ are either {\em cones} with $\ox=0$ or {\em convex} near
$\ox$.

{\bf (iii)} As we show below, the weak contingent extremality
conditions {\em imply} the limiting extremality conditions in any
reflexive space $X$ and also in Asplund spaces under a certain
additional assumption, which is automatic under reflexivity. Thus
the contingent extremality conditions imply the limiting ones in
finite dimensions. The opposite implication does {\em not hold}
even for two sets in $\R^2$. To illustrate it, consider the two
sets from Example~\ref{exam-ext}(i) for which $\ox=(0,0)$ is a
local extremal point in the usual sense, and hence the limiting
extremality conditions hold due to \cite[Theorem~2.8]{m-book1}.
However, it is easy to see that the contingent extremality
conditions are violated for this system.}
\end{Remark}

Observe that for the case of {\em finitely many} sets
$\{\O_1,\ldots,\O_m\}$ the limiting extremality conditions of
Definition~\ref{Def:EPinf}(c) correspond to the generalized Euler
equation in the {\em exact extremal principle} of
\cite[Definition~2.5(iii)]{m-book1} applied to local extremal
points of sets. A natural version of the ``fuzzy" Euler equation
in the {\em approximate extremal principle} of
\cite[Definition~2.5(ii)]{m-book1} for the case of a {\em
countable} set system $\{\O_i\}_{i\in\N}$ at
$\ox\in\cap_{i=1}^\infty\O_i$ can be formulated as follows: for
any $\ve>0$ there are
\begin{equation}\label{aec}
x_i\in\O_i\cap(\ox+\ve\B)\;\mbox{ and }\;x^*_i\in\Hat N(x_i;\O_i)+\frac{1}{2^i}\ve\B^*,\quad i\in\N,
\end{equation}
such that the relationships in \eqref{eq:EPinf} is satisfied. It
turns out that such a countable version of the approximate
extremal principle always {\em holds trivially}, at least in
Asplund spaces, for {\em any} system of closed sets
$\{\O_i\}_{i\in\N}$ at {\em every} boundary point $\ox$ of
infinitely many sets $\O_i$.

\begin{Proposition}\label{aep-tr} {\bf (triviality of the approximate extremality conditions for countable set systems).}  Let
$\{\O_i\}_{i\in\N}$ be a countable system of sets closed around
some point $\ox\in\cap_{i=1}^\infty\O_i$, and let $\ve>0$. Assume
that for infinitely many $i\in\N$ there exist
$x_i\in\O_i\cap(\ox+\ve\B)$ such that $\Hat N(x_i;\O_i)\ne\{0\}$;
this is the case when $X$ is Asplund and $\ox$ belongs to the
boundary of infinitely many sets $\O_i$. Then we always have
$\{x^*_i\}_{i\in\N}$ satisfying conditions \eqref{eq:EPinf} and
\eqref{aec}.
\end{Proposition}
{\bf Proof.} Observe first that the fulfillment of the assumption
made in the proposition for the case of Asplund spaces follows
from the density of Fr\'echet normals on boundaries of closed sets
in such spaces; see, e.g., \cite[Corollary~2.21]{m-book1}. To
proceed further, fix $\ve>0$ and find $j\in\N$ so large that
\begin{equation*}
\disp\frac{\sqrt{2^j}}{2^{j-1}}\le\frac{1}{2}\ve\;\mbox{ and
}\;\Hat N(x_j;\O_j)\ne\{0\}\;\mbox{ with }
\;x_j\in\O_j\cap(\ox+\ve\B).
\end{equation*}
This allows us to get $0\ne x^*_j\in\Hat N(x_j;\O_j)$ such that $\|x^*_j\|=\sqrt{2^j}$ and then choose
\begin{equation*}
\begin{aligned}
&x^*_1:=-\frac{1}{2^{j-1}}x^*_j\in 0+
\frac{1}{2}\ve\B^*\subset\Hat N(x_1;\O_2)+\frac{1}{2}\ve\B^*,\quad
x^*_j\in\Hat N(x_j;\O_j)+\frac{1}{2^j}\ve\B^*,\\
&\mbox{and }\;x^*_i:=0\in\Hat
N(x_i;\O_i)+\frac{1}{2^i}\ve\B^*\;\mbox{ for all }\;i\ne 1,j.
\end{aligned}
\end{equation*}
Thus we have the sequence $\{x^*_i\}_{i\in\N}$ satisfying \eqref{aec} and the relationships
\begin{equation*}
\sum_{i=1}^\infty\frac{1}{2^i}x^*_i=\frac{1}{2}\Big(-\frac{1}{2^{j-1}}x^*_j\Big)+0+\ldots+\frac{1}{2^j}x^*_j
+\ldots=0,\quad\sum_{i=1}^\infty\frac{1}{2^i}\|x^*_i\|^2>1,
\end{equation*}
which give \eqref{eq:EPinf} and complete the proof of the proposition. $\h$

\section{Conic Extremal Principle for Countable Systems of Sets}
\setcounter{equation}{0}

This section addresses the {\em conic extremal principle} for
countable systems of cones in finite-dimensional spaces. This is
the first extremal principle for infinite systems of sets, which
ensures the fulfillment of the conic extremality conditions of
Definition~\ref{Def:EPinf}(a) for a conic extremal system at the
origin under a natural nonoverlapping assumption. We present a
number of examples illustrating the results obtained and the
assumptions made.

To derive the main result of this section, we extend the {\em
method of metric approximations} initiated in \cite{Mor76} to the
case of countable systems of cones; cf.\ an essentially different
realization of this method in the proof of the extremal principle
for local extremal points of finitely many sets in $\R^n$ given in
\cite[Theorem~2.8]{m-book1}. First observe an elementary fact
needed in what follows.

\begin{Lemma}\label{Lem:FderSeries} {\bf (series differentiability).} Let $\|\cdot\|$ be the
usual Euclidian norm in $\R^n$, and let
$\{z_i\}_{i\in\N}\subset\R^n$ be a bounded sequence. Then a
function $\ph\colon\R^n\to\R$ defined by
\begin{equation*}
\ph(x):=\sum_{i=1}^\infty\frac{1}{2^i}\big\|x-z_i\big\|^2,\quad x\in\R^n,
\end{equation*}
is continuously differentiable on $\R^n$ with the derivative
\begin{equation*}
\nabla\ph(x)=\sum_{i=1}^\infty\frac{1}{2^{i-1}}\big(x-z_i\big),\quad x\in\R^n.
\end{equation*}
\end{Lemma}
{\bf Proof}. It is easy to see that both series above converge for
every $x\in\R^n$. Taking further any $u,\xi\in\R^n$ with the norm
$\|\xi\|$ sufficiently small, we have
\begin{equation*}
\|u+\xi\|^2-\|u\|^2-2\la u,\xi\ra=\|u\|^2+2\la u,\xi\ra+\|\xi\|^2-\|u\|^2-2\la u,\xi\ra=\|\xi\|^2=o(\|\xi\|).
\end{equation*}
Thus it follows for any $x\in\R^n$ and $y$ close to $x$ that
\begin{align*}
\ph(y)-\ph(x)-\Big\la\nabla\ph(x),y-x\Big\ra&=
\sum_{i=1}^\infty\frac{1}{2^i}\Big[\|y-z_i\|^2-\|x-z_i\|^2
-2\big\la x-z_i,y-x\big\ra\Big]\\
&=\sum_{i=1}^\infty\frac{1}{2^i}\,\|y-x\|^2=o(\|y-x\|),
\end{align*}
which justifies that $\nabla\ph(x)$ is the derivative of $\ph$ at
$x$, which is obviously continuous on $\R^n$. $\h$\vspace*{0.05in}

Here is the extremal principle for a countable systems of cones,
which plays a crucial role in the subsequent applications of this
paper and its continuation \cite{m-hung}.

\begin{Theorem}\label{Thm:CEP}{\bf (conic extremal principle in finite dimensions).} Let $\{\Lm_i\}_{i\in\N}$
be an extremal
system of closed cones in $X=\R^n$ satisfying the {\sc nonoverlapping condition}
\begin{equation}\label{eq:CEP-QC}
\bigcap_{i=1}^{\infty}\Lm_i=\{0\}.
\end{equation}
Then the conic extremal principle holds, i.e., there are $x^*_i\in N(0;\Lm_i)$ for $i=1,2,\ldots$ such that
\begin{equation*}
\sum_{i=1}^{\infty}\frac{1}{2^i}x^*_i=0\quad\mbox{and}\quad\sum_{i=1}^{\infty}\frac{1}{2^i}\|x^*_i\|^2=1.
\end{equation*}
Moreover, one can find $w_i\in\Lm_i$ for which $x^*_i\in\Hat N(w_i;\Lm_i)$, $i=1,2,\ldots$.
\end{Theorem}
{\bf Proof.} Pick a bounded sequence $\{a_i\}_{i\in\N}\subset\R^n$ from
Definition~\ref{ext-def}(a) satisfying
\begin{equation*}
\bigcap_{i=1}^{\infty}\big(\Lm_i-a_i\big)=\emp
\end{equation*}
and consider the unconstrained optimization problem:
\begin{equation}\label{P1k}
{\rm minimize }\;\ph(x):=\left[\sum_{i=1}^\infty\frac{1}{2^i}\dist^2(x+a_i;\Lm_i)\right]^\frac{1}{2},\quad x\in\R^n.
\end{equation}
Let us prove that problem (\ref{P1k}) has an optimal solution. Since the function $\ph$ in \eqref{P1k} is continuous
on $\R^n$ due the continuity of the distance function and the uniform convergence of the series therein, it suffices
to show that there is $\al>0$ for which the nonempty level set $\{x\in\R^n|\;\ph(x)\le\inf_x\ph+\al\}$ is bounded and
then to apply the classical Weierstrass theorem. Suppose by the contrary that the level sets are unbounded whenever
$\al>0$, for any $k\in\N$ find $x_k\in\R^n$ satisfying
\begin{equation*}
\|x_k\|>k\;\mbox{ and }\;\ph(x_k)\le\inf_x\ph+\frac{1}{k}.
\end{equation*}
Setting $u_k:=x_k/\|x_k\|$ with $\|u_k\|=1$ and taking into account that all $\Lm_i$ are cones, we get
\begin{equation}\label{ser}
\frac{1}{\|x_k\|}\ph(x_k)=\left[\sum_{i=1}^\infty\frac{1}{2^i}\dist^2\Big(u_k+\frac{a_i}{\|x_k\|};\Lm_i\Big)\right]^
\frac{1}{2} \le\frac{1}{\|x_k\|}\Big(\inf_x\ph+\frac{1}{k}\Big)\to
0\;\mbox{ as }\;k\to\infty.
\end{equation}
Furthermore, there is $M>0$ such that for large $k\in\N$
we have
\begin{equation*}
\dist\Big(u_k+\frac{a_i}{\|x_k\|};\Lm_i\Big)\le\Big\|u_k+\frac{a_i}{\|x_k\|}\Big\|\le M.
\end{equation*}
Without relabeling, assume $u_k\to u$ as $k\to\infty$ with some
$u\in\R^n$. Passing now to the limit as $k\to\infty$ in
\eqref{ser} and employing the uniform convergence of the series
therein and the fact that $a_i/\|x_k\|\to 0$ uniformly in $i\in\N$
due the boundedness of $\{a_i\}_{i\in\N}$, we have
\begin{equation*}
\left[\sum_{i=1}^\infty\frac{1}{2^i}\dist^2(u;\Lm_i)\right]^\frac{1}{2}=0.
\end{equation*}
This implies by the closedness of the cones $\Lm_i$ and the {\em
nonoverlapping condition} \eqref{eq:CEP-QC} of the theorem that
$u\in\bigcap_{i=1}^\infty\Lm_i=\{0\}$. The latter is impossible
due to $\|u\|=1$, which contradicts our intermediate assumption on
the unboundedness of the level sets for $\ph$ and thus justifies
the existence of an optimal solution $\Tilde x$ to problem
\eqref{P1k}.

Since the system of closed cones $\{\Lm_i\}_{i\in\N}$ is {\em
extremal at the origin}, it follows from the construction of $\ph$
in \eqref{P1k} that $\ph(\Tilde x)>0$. Taking into account the
nonemptiness of the projection $\Pi(x;\Lm)$ of $x\in\R^n$ onto an
arbitrary closed set $\Lm\subset\R^n$, pick any $w_i\in\Pi(\Tilde
x+a_i;\Lm_i)$ as $i\in\N$ and observe from
Proposition~\ref{Prop:N(Lm)} above and the proof of
\cite[Theorem~1.6]{m-book1} that
\begin{equation}\label{pr1}
\Tilde x+a_i-w_i\in \Pi^{-1}(w_i;\Lm_i)-w_i\subset\Hat N(w_i;\Lm_i)\subset N(0;\Lm_i).
\end{equation}
Furthermore, the sequence $\{a_i-w_i\}_{i\in\N}$ is bounded in $\R^n$ due to
\begin{equation*}
\|x+a_i-w_i\|=\dist(x+a_i;\Lm_i)\le\|x+a_i\|.
\end{equation*}

Next we consider another unconstrained optimization problem:
\begin{equation}
\label{P2k}{\rm minimize }\
\psi(x):=\left[\sum_{i=1}^\infty\frac{1}{2^i}\|x+a_i-w_i\|^2\right]^\frac{1}{2},\quad x\in\R^n.
\end{equation}
It follows from $\psi(x)\ge\ph(x)\ge\ph(\Tilde x)=\psi(\Tilde x)$
for all $x\in\R^n$ that problem \eqref{P2k} has the same optimal
solution $\Tilde x$ as \eqref{P1k}. The main difference between
these two problems is that the cost function $\psi$ in \eqref{P2k}
is {\em smooth} around $\Tilde x$ by Lemma~\ref{Lem:FderSeries},
the smoothness of the function $\sqrt{t}$ around nonzero points,
and the fact that $\psi(\Tilde x)\ne 0$ due to the cone
extremality.  Applying now the {\em classical Fermat rule} to the
smooth unconstrained minimization problem (\ref{P2k}) and using
the derivative calculation in Lemma~\ref{Lem:FderSeries}, we
arrive at the relationships
\begin{equation}\label{pr2}
\nabla\psi(\Tilde x)=\sum_{i=1}^\infty\frac{1}{2^i}x^*_i=0\;\mbox{
with }\;x^*_i:=\frac{1}{\psi(\Tilde x)}\Big(\Tilde x+a_i-w_i\Big),
\quad i\in\N.
\end{equation}
The latter implies by \eqref{pr1} that $x^*_i\in\Hat N(w_i;\Lm_i)\subset N(0;\Lm_i)$ for all $i\in\N$. Furthermore,
it follows from the constructions of $x^*_i$ in \eqref{pr2} and of $\psi$ in \eqref{P2k} that
\begin{equation*}
\sum_{i=1}^\infty\frac{1}{2^i}\|x^*_i\|^2=1,
\end{equation*}
which thus completes the proof of the theorem. $\h$\vspace*{0.05in}

In the remaining part of this section, we present three examples
showing that all the assumptions made in Theorem~\ref{Thm:CEP}
(nonoverlapping, finite dimension, and conic structure) are {\em
essential} for the validity of this result.

\begin{Example}\label{Ex:fail-EP1} {\bf (nonoverlapping condition is essential).} {\rm Let us show that the conic
extremal principle may fail for countable systems of {\em convex}
cones in $\R^2$ if the nonoverlapping condition \eqref{eq:CEP-QC}
is violated. Define the convex cones $\Lm_i\subset\R^2$ as
$i\in\N$ by
\begin{equation*}
\Lm_1:=\R\times \R_+\;\mbox{ and }\;\Lm_i:=\big\{(x,y)\in\R^2\big|\;y\le\frac{x}{i}\big\}\;\mbox{ for }\;i=2,3,\ldots.
\end{equation*}
Observe that for any $\nu>0$ we have
\begin{equation*}
\Big(\Lm_1+(0,\nu)\Big)\bigcap_{k=2}^\infty\Lm_k=\emp,
\end{equation*}
which means that the cone system $\{\Lm_i\}_{i\in\N}$ is extremal at the origin. On the other hand,
\begin{equation*}
\bigcap_{i=1}^\infty\Lm_i=\R_+\times\{0\},
\end{equation*}
i.e., the nonoverlapping condition \eqref{eq:CEP-QC} is violated.
Furthermore, we can easily compute the corresponding normal cones
by
\begin{equation*}
N(0;\Lm_1)=\big\{\lm(0,-1)\big|\;\lm\ge 0\big\}\;\mbox{ and }\;
N(0;\Lm_i)=\big\{\lm(-1,i)\big|\;\lm\ge 0\big\},\;\;i=2,3,\ldots.
\end{equation*}
Taking now any $x^*_i\in N(0;\Lm_i)$ as $i\in\N$, observe the equivalence
\begin{equation*}
\Big[\sum_{i=1}^{\infty}\frac{1}{2^i}x^*_i=0\Big]\Longleftrightarrow\Big[\frac{\lm_1}{2}\big(0,-1\big)+
\sum_{i=2}^{\infty}\frac{\lm_i}{2^i}\big(-1,i\big)=0\;\mbox{ with }\;\lm_i\ge 0\;\mbox{ as }\;i\in\N\Big].
\end{equation*}
The latter implies that $\lm_i=0$ and hence $x^*_i=0$ for all
$i\in\N$. Thus the nontriviality condition in \eqref{eq:EPinf} is
not satisfied, which shows that the conic extremal principle fails
for this system.}
\end{Example}

\begin{Example}\label{ex2} {\bf (conic structure is essential).} {\rm If all the sets $\O_i$ for $i\in\N$ are
{\em convex} but some of them are {\em not cones}, then the
equivalent extremality conditions of
Definition~\ref{Def:EPinf}(b,c) are natural extensions of the
conic extremality conditions in Theorem~\ref{Thm:CEP}. We show
nevertheless that the corresponding extension of the conic
extremal principle under the nonoverlapping requirement
\begin{equation}\label{overlap}
\bigcap_{i=1}^\infty\O_i=\{0\}
\end{equation}
{\em fails} without imposing a conic structure on all the sets
involved. Indeed, consider a countable system of closed and convex
sets in $\R^2$ defined by
\begin{equation*}
\O_1:=\big\{(x,y)\in\R^2\big|\;y\ge x^2\big\}\;\mbox{ and }\;\O_i:=\big\{(x,y)\in\R^2\big|\;y\le\frac{x}{i}\big\}
\;\mbox{ for }\;i=2,3,\ldots.
\end{equation*}
We can see that only the set $\O_1$ is not a cone and that the
nonoverlapping requirement \eqref{overlap} is satisfied.
Furthermore, the system $\{\O_i\}_{i\in\N}$ is extremal at the
origin in the sense that \eqref{cone-emp} holds. However, the
arguments similar to Example~\ref{Ex:fail-EP1} show that the
extremality conditions \eqref{eq:EPinf} with $x^*_i\in N(0;\O_i)$
as $i\in\N$ fail to fulfill. Note that, as shown in Section~7,
both contingent and limiting extremal principles hold for
countable systems of general nonconvex sets if nonoverlapping
condition \eqref{overlap} is replaced by another one reflecting
the {\em contingent extremality}.}
\end{Example}

\begin{Example}\label{ex3}{\bf (failure of the conic extremal principle in infinite dimensions).} {\rm The last example
demonstrates that the conic extremal principle of
Theorem~\ref{Thm:CEP} with the nonoverlapping condition
\eqref{eq:CEP-QC} may fail for countable systems of {\em convex
cones} (in fact, half-spaces) in an arbitrary infinite-dimensional
{\em Hilbert space}. To proceed, consider a Hilbert space $X$ with
the orthonormal basis $\{e_i|\;i\in\N\}$ and define a countable
system of closed half-spaces by
\begin{equation*}
\Lm_1:=\big\{x\in X\big|\;\la x,e_1\ra\le 0\big\}\;\mbox{ and }\;\Lm_i:=\big\{x\in X\big|\;\la x,
e_{i}-e_{i-1}\ra\le 0\big\}\;\mbox{ for }\;i=2,3,\ldots.
\end{equation*}
It is easy to compute the corresponding normal cones to the above sets:
\begin{equation*}
N(0;\Lm_1)=\big\{\lm e_1\big|\;\lm\ge 0\big\}\;\mbox{ and
}\;N(0;\Lm_i)=\big\{\lm(e_{i}-e_{i-1})\big|\;\lm\ge 0\big\}\;
\mbox{ for }\;i=2,3,\ldots.
\end{equation*}
Now let us check that the nonoverlapping condition \eqref{eq:CEP-QC} is satisfied. Indeed, picking any point
\begin{equation*}
x=\sum_{i=1}^\infty\al_ie_i\in\bigcap_{i=1}^\infty\Lm_i,
\end{equation*}
we have $\al_1=\la x,e_1\ra\le 0$ and $\al_i=\la x,e_{i}\ra\le\la
x,e_{i-1}\ra=\al_{i-1}$ for $i=2,3,\ldots$. This clearly leads to
$\al_i=0$ for all $i\in\N$, which yields $x=0$ and thus justifies
(\ref{eq:CEP-QC}). The same arguments show that
\begin{equation*}
(\Lm_1-e_1)\cap\bigcap_{i=2}^\infty\Lm_i=\emp,
\end{equation*}
i.e., $\{\Lm_i\}_{i\in\N}$ is a {\em conic extremal system}.
However, the conic extremality conditions of
Definition~\ref{Def:EPinf}(a) {\em fail} for this system. To check
this, suppose that there exist $x^*_i\in N(0;\Lm_i)$ as $i\in\N$
satisfying the relationships
\begin{equation}\label{ex3-1}
\displaystyle\sum_{i=1}^\infty
x^*_i=0\;\mbox{ and }\;\displaystyle\sum_{i=1}^\infty\|x^*_i\|>0.
\end{equation}
By the above structure of $N(0;\Lm_i)$ we have $x^*_1=\lm_1 e_1$
and $x^*_i=\lm_i(e_{i}-e_{i-1})$ as $i=2,3,\ldots$ for some
$\lm_i\ge 0$ as $i\in\N$. Thus the first condition in
\eqref{ex3-1} reduces to
\begin{equation*}
\lm_1e_1+\sum_{i=2}^\infty\lm_i\big(e_{i}-e_{i-1}\big)=0.
\end{equation*}
The latter is possible if either (a): $\lm_i=1$ for all $i\in\N$
or (b): $\lm_i=0$ for all $i\in\N$. Case (a) surely contradicts
the convergence of the series in the second condition of
\eqref{ex3-1} while in case (b) the latter series converges to
zero. Hence the conic extremal principle of Theorem~\ref{Thm:CEP}
does not hold in this infinite-dimensional setting.}
\end{Example}

\section{Fr\'echet Normals to Countable Intersections of Cones}
\setcounter{equation}{0}

In this section we present applications of the conic extremal
principle established in Theorem~\ref{Thm:CEP} to deriving several
representations, under appropriate assumptions, of Fr\'echet
normals to {\em countable intersections} of cones in
finite-dimensional spaces. These calculus results are certainly of
their independent interest while their are largely employed in
\cite{m-hung} to problems of semi-infinite programming and
multiobjective optimization.

To begin with, we introduce the following qualification condition
for countable systems of cones formulated in terms of limiting
normals \eqref{nc}, which plays a significant role in deriving the
results of this section as well as in the subsequent applications
given in \cite{m-hung}.

\begin{Definition}\label{eq:QC-N(Lm)}{\bf (normal qualification condition for countable systems of cones).}
Let $\{\Lm_i\}_{i\in\N}$ be a countable system of closed cones in
$X$. We say that it satisfies the {\sc normal qualification
condition} at the origin if
\begin{equation}\label{nqc}
\Big[\sum_{i=1}^{\infty}x^*_i=0,\;\;x^*_i\in
N(0;\Lm_i)\Big]\Longrightarrow\big[x^*_i=0,\;\;i\in\N\big].
\end{equation}
\end{Definition}

This definition corresponds to the normal qualification condition
of \cite{m-book1} for finite systems of sets; seethe discussions
and various applications of the latter condition therein. We refer
the reader to \cite{m-hung} for a nonconic version of \eqref{nqc},
its relationships with other qualification conditions for
countable systems of sets, and sufficient conditions for its
validity that equally apply to both conic and nonconic versions.
In this section we use the normal qualification condition of
Definition~\ref{eq:QC-N(Lm)} to represent Fr\'echet normals to
countable intersections of cones in terms of limiting normals to
each of the sets involved. Let us start with the following
``fuzzy" intersection rule at the origin.

\begin{Theorem}\label{Thm:Fuzzy}{\bf (fuzzy intersection rule for \F\ normals to countable intersections of cones).}
Let $\{\Lm_i\}_{i\in\N}$ be a countable system of arbitrary closed
cones in $\R^n$ satisfying the normal qualification condition
\eqref{nqc}. Then given a Fr\'echet normal $x^*\in\Hat
N\big(0;\bigcap_{i=1}^\infty\Lm_i\big)$ and a number $\ve>0$,
there are limiting normals $x^*_i\in N(0;\Lm_i)$ as $i\in\N$ such
that
\begin{equation}\label{fuz}
x^*\in\sum_{i=1}^\infty \frac{1}{2^i}x^*_i+\ve\B^*.
\end{equation}
\end{Theorem}
{\bf Proof}. Fix $x^*\in\Hat N\Big(0;\bigcap_{i=1}^\infty\Lm_i\Big)$
and $\ve>0$. By definition \eqref{eq:e-nor} of \F\ normals we have
\begin{equation}\label{fuz1}
\la x^*,x\ra-\ve\|x\|<0\;\mbox{ whenever }\;x\in\bigcap_{i=1}^\infty\Lm_i\setminus\{0\}.
\end{equation}
Define a countable system of closed cones in $\R^{n+1}$ by
\begin{equation}\label{o1}
O_1:=\big\{(x,\al)\in\R^n\times\R\big|\;x\in\Lm_1,\;\al\le\la x^*,x\ra-\ve\|x\|\big\}\;\mbox{ and }\;
O_i:=\Lm_i\times\R_+\;\mbox{ for }\;i=2,3,\ldots.
\end{equation}
Let us check that all the assumptions for the validity of the
conic extremal principle in Theorem~\ref{Thm:CEP} are satisfied
for the system $\{O_i\}_{i\in\N}$. Picking any
$(x,\al)\in\bigcap_{i=1}^\infty O_i$, we have
$x\in\bigcap_{i=1}^\infty\Lm_i$ and $\al\ge 0$ from the
construction of $\O_i$ as $i\ge 2$. This implies in fact that
$(x,\al)=(0,0)$. Indeed, supposing $x\ne 0$ gives us by
\eqref{fuz1} that
\begin{equation*}
0\le\al\le\la x^*,x\ra-\ve\|x\|<0,
\end{equation*}
which is a contradiction. On the other hand, the inclusion
$(0,\al)\in O_1$ yields that $\al\le 0$ by the construction of
$O_1$, i.e., $\al=0$. Thus the {\em nonoverlapping condition}
\begin{equation*}
\bigcap_{i=1}^\infty O_i=\{(0,0)\}
\end{equation*}
holds for $\{O_i\}_{i\in\N}$. Similarly we check that
\begin{equation}\label{eq:Fz1}
\Big(O_1-(0,\gg)\Big)\cap\bigcap_{i=2}^\infty O_i=\emp\;\mbox{ for any fixed }\;\gg>0,
\end{equation}
i.e., $\{O_i\}_{i\in\N}$ is a {\em conic extremal system} at the
origin. Indeed, violating \eqref{eq:Fz1} means he existence of
$(x,\al)\in\R^n\times\R$ such that
\begin{equation*}
(x,\al)\in\Big[O_1-(0,\gg)\Big]\cap\bigcap_{i=2}^\infty O_i,
\end{equation*}
which implies that $x\in\bigcap_{i=1}^\infty O_i$ and $\al\ge 0$. Then by the
construction of $O_1$ in \eqref{o1} we get
\begin{equation*}
\gg+\al\le\la x^*,x\ra-\ve\|x\|\le 0,
\end{equation*}
a contradiction due the positivity of $\gg$ in \eqref{eq:Fz1}.

Applying now the second conclusion of Theorem~\ref{Thm:CEP} to the
system $\{O_i\}_{i\in\N}$ gives us the pairs $(w_i,\al_i)\in O_i$
and $(x^*_i,\lm_i)\in\Hat N\big((w_i,\al_i);O_i\big)$ as $i\in\N$
satisfying the relationships
\begin{equation}\label{fuz2}
\sum_{i=1}^{\infty}\frac{1}{2^i}\big(x^*_i,\lm_i\big)=0\quad\mbox{and}\quad\sum_{i=1}^{\infty}\frac{1}{2^i}\big\|
(x^*_i,\lm_i)\big\|^2=1.
\end{equation}
It immediately follows from the constructions of $O_i$ as $i\ge 2$ in \eqref{o1} that
$\lm_i\le 0$ and $x^*_i\in\Hat N(w_i;\Lm_i)$; thus $x^*_i\in N(0;\Lm_i)$ for $i=2,3,\ldots$ by
Proposition~\ref{Prop:N(Lm)}. Furthermore,
we get
\begin{equation}\label{eq:FzLimsup}
\limsup_{(x,\al)\st{O_1}{\to}(w_1,\al_1)}\frac{\la
x^*_1,x-w_1\ra+\lm_1(\al-\al_1)}{\|x-w_1\|+|\al-\al_1|}\le 0
\end{equation}
by the definition of \F\ normals to $O_1$ at $(w_1,\al_1)\in O_1$ with
$\lm_1\ge 0$ and
\begin{equation}\label{eq:Fz2}
\al_1\le\la x^*,w_1\ra -\ve\|w_1\|
\end{equation}
by the construction of $O_1$. Examine next the two possible cases in \eqref{fuz2}: $\lm_1=0$ and $\lm_1>0$.\\[1.1ex]
{\bf Case~1}: $\lm_1=0$. If inequality (\ref{eq:Fz2}) is strict in this case, find a neighborhood $U$ of $w_1$ such
that
\begin{equation*}
\al_1<\la x^*,x\ra-\ve\|x\|\  \mbox{ for all }\ x\in U,
\end{equation*}
which ensures that $(x,\al_1)\in O_1$ for all $x\in \Lm_1\cap U$. Substituting $(x,\al_1)$ into (\ref{eq:FzLimsup})
gives us
\begin{equation*}
\limsup_{x\st{\Lm_1}{\to}w_1}\frac{\la x^*_1,x-w_1\ra}{\|x-w_1\|}\le 0,
\end{equation*}
which means that $x^*_1\in\Hat N(w_1;\Lm_1)$. If (\ref{eq:Fz2}) holds as equality, we put $\al:=\la x^*,x\ra -\ve\|x\|$
and get
\begin{equation*}
|\al-\al_1|=\big|\la x^*,x-w_1\ra +\ve(\|w_1\|-\|x\|)\big|\le\big(\|x^*\|+\ve\big)\|x-w_1\|.
\end{equation*}
Furthermore, it follows from (\ref{eq:FzLimsup}) that
\begin{equation*}
\limsup_{(x,\al)\st{O_1}{\to}(w_1,\al_1)}
\frac{\la x^*_1,x-w_1\ra}{\|x-w_1\|+|\al-\al_1|}\le 0.
\end{equation*}
Thus for any $\nu>0$ sufficiently small and $\al$ chosen above, we have
\begin{equation*}
\la x^*_1,x-w_1\ra\le\nu\big(\|x-w_1\|+|\al-\al_1|\big)\le\nu\big(1+\|x^*\|+\ve\big)\|x-w_1\|
\end{equation*}
whenever $x\in\Lm_1$ is sufficiently closed to $w_1$. The latter yields that
\begin{equation*}
\disp\limsup_{x\st{\Lm_1}{\to}w_1}\frac{\la
x^*_1,x-w_1\ra}{\|x-w_1\|}\le 0,\;\mbox{  i.e., }\;x^*_1\in\Hat N(w_1;\Lm_1).
\end{equation*}
Thus in both cases of the strict inequality and equality in (\ref{eq:Fz2}),
we justify that $x^*_1\in\Hat N(w_1;\Lm_1)$ and thus $x^*_1\in N(0;\Lm_1)$ by
Proposition~\ref{Prop:N(Lm)}. Summarizing the above discussions gives us
\begin{equation*}
x^*_i\in N(0;\Lm_i)\;\mbox{ and }\;\lm_i=0\;\mbox{ for all }\;i\in\N
\end{equation*}
in Case~1 under consideration. Hence it follows from \eqref{fuz2} that there are $\Tilde x^*_i:=(1/2^i)x^*_i\in
N(0;\Lm_i)$ as $i\in\N$, not equal to zero simultaneously, satisfying
\begin{equation*}
\sum_{i=1}^\infty\Tilde x^*_i=0.
\end{equation*}
This contradicts the normal qualification condition \eqref{nqc} and thus shows that the case of $\lm_1=0$ is
actually {\em not possible} in \eqref{eq:Fz2}.
\\[1.1ex]
{\bf Case~2:} $\lm_1>0$. If inequality (\ref{eq:Fz2}) is strict, put $x=w_1$ in (\ref{eq:FzLimsup}) and get
\begin{equation*}
\limsup_{\al\to\al_1}\frac{\lm_1(\al-\al_1)}{|\al-\al_1|}\le 0.
\end{equation*}
That yields $\lm_1=0$, a contradiction. Hence it remains to consider the case when (\ref{eq:Fz2}) holds as equality.
To proceed, take $(x,\al)\in O_1$ satisfying
\begin{equation*}
x\in\Lm_1\setminus\{w_1\}\;\mbox{ and }\;\al=\la x^*,x\ra-\ve\|x\|.
\end{equation*}
By the equality in (\ref{eq:Fz2}) we have
\begin{equation*}
\al-\al_1=\la x^*,x-w_1\ra+\ve(\|w_1\|-\|x\|)\;\mbox{ and thus }\;
|\al-\al_1|\leq(\|x^*\|+\ve)\|x-w_1\|.
\end{equation*}
On the other hand, it follows from (\ref{eq:FzLimsup}) that for any
$\gg>0$ sufficiently small there exists a neighborhood $V$ of $w_1$ such that
\begin{equation}\label{fuz3}
\la x^*_1,x-w_1\ra+\lm_1(\al-\al_1)\le\lm_1\gg\ve\big(\|x-w_1\|+|\al-\al_1|\big)
\end{equation}
whenever $x\in\Lm_1\cap V$. Substituting $(x,\al)$ with $x\in\Lm_1\cap V$ into \eqref{fuz3} gives us
\begin{align*}
\la x^*_1,x-w_1\ra +\lm_1(\al-\al_1)
&=\la x^*_1+\lm_1 x^*,x-w_1\ra +\lm_1\ve(\|w_1\|-\|x\|)\\
&\le\lm_1\gg\ve(\|x-w_1\|+|\al-\al_1|)\\
&\le\lm_1\gg\ve\big[\|x-w_1\|+(\|x^*\|+\ve)\|x-w_1\|\big]\\
&=\lm_1\gg\ve\big(1+\|x^*\|+\ve\big)\|x-w_1\|.
\end{align*}
It follows from the above that for small $\gg>0$ we have
\begin{equation*}
\la x^*_1+\lm_1 x^*,x-w_1\ra +\lm_1\ve(\|w_1\|-\|x\|)
\le\lm_1\ve\|x-w_1\|
\end{equation*}
and thus arrive at the estimates
\begin{equation*}
\la x^*_1+\lm_1x^*,x-w_1\ra\le\lm_1\ve\|x-w_1\|+\lm_1\ve(\|x\|-\|w_1\|)\le 2\lm_1\ve\|x-w_1\|
\end{equation*}
for all $x\in\Lm_1\cap V$. The latter implies by definition \eqref{eq:e-nor} of $\ve$-normals that
\begin{equation}\label{fuz4}
x^*_1+\lm_1x^*\in\Hat N_{2\lm_1\ve}(w_1;\Lm_1).
\end{equation}
Furthermore, it is easy to observe from the above choice of
$\lm_1$ and the structure of $O_1$ in \eqref{o1} that $\lm_1\le
2+2\ve$. Employing now the representation of $\ve$-normals in
\eqref{fuz4} from \cite[formula (2.51)]{m-book1} held in finite
dimensions, we find $v\in\Lm_1\cap(w_1+2\lm_1\ve\B)$ such that
\begin{equation}\label{fuz5}
x^*_1+\lm_1x^*\in\Hat N(v;\Lm_1)+2\lm_1\ve\B^*\subset N(0;\Lm_1)+2\lm_1\ve\B^*.
\end{equation}
Since $\lm_1>0$ in the case under consideration and by
$\displaystyle-x^*_1=2\sum_{i=2}^\infty\frac{1}{2^i}x^*_i$ due to
the first equality in \eqref{fuz2}, it follows from \eqref{fuz5}
that
\begin{equation*}
x^*\in N(0;\Lm_1)+\frac{2}{\lm_1}\sum_{i=2}^\infty\frac{1}{2^i}x^*_i+2\ve\B^*,
\end{equation*}
and hence there exists $\Tilde x^*_1\in N(0;\Lm_1)$ such that
\begin{equation*}
x^*\in\sum_{i=1}^\infty\frac{1}{2^i}\Tilde x^*_i+2\ve\B^*\;\mbox{
with }\;\Tilde x^*_i:=\frac{2x^*_i}{\lm_1}\in N(0;\Lm_i)\; \mbox{
for }\;i=2,3,\ldots.
\end{equation*}
This justifies \eqref{fuz} and completes the proof of the theorem. $\h$\vspace*{0.1in}

Our next result shows that we can put $\ve=0$ in representation
\eqref{fuz} under an additional assumption on Fr\'echet normals to
cone intersections.

\begin{Theorem}\label{Thm:(r)Fuzzy} {\bf (refined representation of Fr\'echet normals to countable intersections of
cones).} Let $\{\Lm_i\}_{i\in\N}$ be a countable system of
arbitrary closed cones in $\R^n$ satisfying the normal
qualification condition \eqref{nqc}. Then for any Fr\'echet normal
$x^*\in\Hat N\Big(0;\bigcap_{i=1}^\infty\Lm_i\Big)$ satisfying
\begin{equation}\label{i1}
\la x^*,x\ra<0\;\mbox{ whenever }\;x\in\bigcap_{i=1}^\infty\Lm_i\setminus\{0\}
\end{equation}
there are limiting normals $x^*_i\in N(0;\Lm_i)$, $i=1,2,\ldots$, such that
\begin{equation}\label{i2}
x^*=\sum_{i=1}^\infty\frac{1}{2^i}x^*_i.
\end{equation}
\end{Theorem}
{\bf Proof.} Fix a Fr\'echet normal $x^*\in\Hat
N\Big(0;\bigcap_{i=1}^\infty\Lm_i\Big)$ satisfying condition
\eqref{i1} and construct a countable system of closed cones in
$\R^n\times\R$ by
\begin{equation}\label{o2}
O_1:=\big\{(x,\al)\in\R^n\times\R\big|\;x\in\Lm_1,\;\al\le\la x^*,x\ra\big\}\;\mbox{ and }\;
O_i:=\Lm_i\times\R_+\mbox{ for}\;i=2,3,\ldots.
\end{equation}
Similarly to the proof Theorem~\ref{Thm:Fuzzy} with taking
\eqref{i1} into account, we can verify that all the assumptions of
Theorem~\ref{Thm:CEP} hold. Applying the conic extremal principle
from this theorem gives us pairs $(w_i,\al_i)\in O_i$ and
$(x^*_i,\lm_i)\in\Hat N\big((w_i,\al_i);O_i\big)$ such that the
extremality conditions in \eqref{fuz2} are satisfied. We obviously
get $\lm_i\le 0$ and $x^*_i\in\Hat N(w_i;\Lm_i)$ for
$i=1,2,\ldots$, which ensures that $x^*_i\in N(0;\Lm_i)$ as $i\ge
2$ by Proposition~\ref{Prop:N(Lm)}. It follows furthermore that
for $i=1$ the limiting inequality \eqref{eq:FzLimsup} holds. The
latter implies by the structure of the set $O_1$ in \eqref{o2}
that
\begin{equation}\label{i3}
\lm_1\ge 0\;\mbox{ and }\;\al_1\le\la x^*,w_1\ra.
\end{equation}
Similarly to the proof of Theorem~\ref{Thm:Fuzzy} we consider the
two possible cases $\lm_1=0$ and $\lm_1>0$ in \eqref{i3} and show
that the first case contradicts the normal qualification condition
\eqref{nqc}. In the second case we arrive at representation
\eqref{i2} based on the extremality conditions in \eqref{fuz2} and
the structures of the sets $O_i$ in \eqref{o2}.
$\h$\vspace*{0.05in}

The next theorem in this section provides constructive upper estimates of the Fr\'echet normal cone to countable
intersections of closed cones
in finite dimensions and of its interior via limiting normals to the sets involved at the origin.

\begin{Theorem}\label{Thm:FnorLm} {\bf (Fr\'echet normal cone to countable intersections).} Let $\{\Lm_i\}_{i\in\N}$
be a countable system of arbitrary closed cones in $\R^n$
satisfying the normal qualification condition \eqref{nqc}, and let
$\Lm:=\bigcap_{i=1}^\infty\Lm_i$. Then we have the inclusions
\begin{equation}\label{int}
\inter\Hat N(0;\Lm)\subset\Big\{\sum_{i=1}^{\infty}x^*_i\Big|\;x^*_i\in N(0;\Lm_i)\Big\},
\end{equation}
\begin{equation}\label{cl}
\Hat N(0;\Lm)\subset\cl\Big\{\sum_{i\in I}x^*_i\Big|\;x^*_i\in N(0;\Lm_i),\;I\in{\cal L}\Big\},
\end{equation}
where ${\cal L}$ stands for the collection of all finite subsets of the natural series $\N$.
\end{Theorem}
{\bf Proof.} First we justify inclusion \eqref{int} assuming
without loss of generality that $\inter N(0;\Lm)\ne\emp$. Pick any
$x^*\in\inter\Hat N(0;\Lm)$ and also $\gg>0$ such that
$x^*+3\gg\B^*\subset\Hat N(0;\Lm)$. Then for any
$x\in\Lm\setminus\{0\}$ find $z^*\in\R^n$ satisfying the
relationships
\begin{equation*}
\|z^*\|=2\gg\;\mbox{ and }\;\la z^*,x\ra<-\gg\|x\|.
\end{equation*}
Since $x^*-z^*\in x^*+3\gg\B^*\subset\Hat N(0;\Lm)$, we have $\la x^*-z^*,x\ra\le 0$ and hence
\begin{equation*}
\la x^*,x\ra=\la x^*-z^*,x\ra+\la z^*,x\ra<-\gg\|x\|<0.
\end{equation*}
This allows us to employ Theorem~\ref{Thm:(r)Fuzzy} and thus justify the first inclusion \eqref{int}.

To prove the remaining inclusion \eqref{cl}, pick pick $x^*\in\Hat N(0;\Lm)$ and for any fixed $\ve>0$ apply
Theorem~\ref{Thm:Fuzzy}.
In this way we find $x^*_i\in N(0;\Lm_i)$, $i\in\N$, such that
\begin{equation*}
x^*\in\sum_{i=1}^\infty\frac{1}{2^i}x^*_i+\ve\B^*.
\end{equation*}
Since $\ve>0$ was chosen arbitrarily, it follows that
\begin{equation*}
x^*\in A:=\cl\left\{\sum_{i=1}^\infty\frac{1}{2^i}x^*_i\Big|\;x^*_i\in N(0;\Lm_i)\right\}.
\end{equation*}
Let us finally justify the inclusion
\begin{equation*}\label{C}
A\subset\cl C\;\mbox{ with }\;C:=\Big\{\sum_{i\in I}x^*_i\Big|\;x^*_i\in N(0;\Lm_i),\;I\in{\cal L}\Big\}.
\end{equation*}
To proceed, pick $z^*\in A$ and for any fixed $\ve>0$ find $x^*_i\in N(0;\Lm_i)$ satisfying
\begin{equation*}
\left\|z^*-\sum_{i=1}^\infty\frac{1}{2^i}x^*_i\right\|\le\frac{\ve}{2}.
\end{equation*}
Then choose a number $k\in\N$ so large that
\begin{equation*}
\left\|z^*-\sum_{i=1}^k\frac{1}{2^i}x^*_i\right\|\le\ve.
\end{equation*}
Since $\disp\sum_{i=1}^k\frac{1}{2^i}x^*_i\in C$, we get $(z^*+\ve\B^*)\cap C\ne\emp$, which means that
$z^*\in\cl C$. This justifies \eqref{cl} and completes the proof of the theorem.$\h$\vspace*{0.05in}

Finally in this section, we present a consequence of Theorem~\ref{Thm:FnorLm}, which gives an exact computation
of Fr\'echet normals to countable intersections of cones normally regular at the origin.

\begin{Corollary}\label{n-reg}{\bf (countable intersections of normally regular cones).} In addition to the assumptions
of Theorem~{\rm\ref{Thm:FnorLm}}, suppose that all the cones
$\Lm_i$, $i\in\N$, are normally regular at the origin. Then the
Fr\'echet normal cone to the intersection
$\Lm=\bigcap_{i=1}^\infty\Lm_i$ is computed by
\begin{equation}\label{f-eq}
\Hat N(0;\Lm)=\cl\Big\{\sum_{i\in I}x^*_i\Big|\;x^*_i\in N(0;\Lm_i),\;I\in{\cal L}\Big\}.
\end{equation}
\end{Corollary}
{\bf Proof.} It is easy to check that
\begin{equation*}
\cl\Big\{\sum_{i\in I}x^*_i\Big|\;x^*_i\in\Hat N(0;\Lm_i),\;I\in{\cal L}\Big\}\subset\Hat N(0;\Lm)
\end{equation*}
for arbitrary set systems. Combining this with inclusion
\eqref{cl} of Theorem~\ref{Thm:FnorLm} and the normal regularity
of each cone $\Lm_i$ as $i\in\N$, gives us equality \eqref{f-eq}.
$\h$

\section{Tangential Normal Enclosedness and Approximate Normality}
\setcounter{equation}{0}

In this section we introduce and study two important properties of
tangents cones that are of their own interest while allow us make
a bridge between the extremal principles for cones and the
limiting extremality conditions for arbitrary closed sets at their
tangential extremal points. The main attention is paid to the
contingent and weak contingent cones, which are proved to enjoy
these properties under natural assumptions.

Let us start with introducing a new property of sets that is
formulated in terms of the limiting normal cone \eqref{nc} and
plays a crucial role of what follows.

\begin{Definition}\label{Def:TNE} {\bf (tangential normal enclosedness).} Given a nonempty subset $\O\subset X$ and a
subcone $\Lm\subset X$ of a Banach space $X$, we say that $\Lm$ is
{\sc tangentially normally enclosed} {\rm(TNE)} into $\O$ at a
point $\ox\in\O$ if
\begin{equation}\label{tne}
N(0;\Lm)\subset N(\ox;\O).
\end{equation}
\end{Definition}
The word ``tangential" in Definition~\ref{Def:TNE} reflects the
fact that this normal enclosedness property is applied to
tangential approximations of sets at reference points. Observe
that if the set $\O$ is convex near $\ox$, then its classical
tangent cone at $\ox$ enjoys the TNE property; indeed, in this
case inclusion \eqref{tne} holds as equality. We establish below a
remarkable fact on the validity of the TNE property for the weak
contingent cone \eqref{wbs} to any closed subset of a reflexive
Banach space.\vspace*{0.05in}

To study this and related properties, fix $\O\subset X$ with
$\ox\in\O$ and denote by $\Lm_w:=T_w(\ox;\O)$ the weak contingent
cone to $\O$  at $\ox$ without indicating $\O$ and $\ox$ for
brevity. Given a direction $d\in\Lm_w$, let $\mT^w_d$ be the
collection of all sequences $\{x_k\}\subset\O$ such that
\begin{equation*}
\frac{x_k-\ox}{t_k}\st{w}{\lto}d\;\mbox{ for some }\;t_k\dn 0.
\end{equation*}
It follows from definition \eqref{wbs} of $\Lm_w=T(\ox;\O)$ that $\mT^w_d\ne\emp$ whenever $d\in\Lm_w$.

\begin{Definition}\label{Def:(AN)}{\bf(tangential approximate normality).} We say that $\O\subset X$ has the
{\sc tangential  approximate normality} {\rm(TAN)} property at
$\ox\in\O$ if whenever $d\in\Lm_w$ and $x^*\in\Hat N(d;\Lm_w)$ are
chosen there is a sequence $\{x_k\}\in\mT^w_d$ along which the
following holds: for any $\ve>0$ there exists $\dd\in(0,\ve)$ such
that
\begin{equation}\label{eq:(AN)}
\limsup_{k\to\infty}\Big[\sup\Big\{\frac{\la
x^*,z-x_k\ra}{t_k}\Big|\;z\in\O\cap (x_k+t_k\dd\B)\Big\}\Big] \le
2\ve\dd,
\end{equation}
where $t_k\dn0$ is taken from the construction of $\mT^w_d$.
\end{Definition}

The meaning of this property that gives the name is as follows:
any $x^*\in\Hat N(d;\Lm_w)$ for the tangential approximation of
$\O$ at $\ox$ behaves approximately like a true normal at
appropriate points $x_k$ near $\ox$. It occurs that the TAN
property holds for any closed subset of a reflexive Banach space.
The next proposition provides even a stronger result.

\begin{Proposition}\label{Lem:(AN)} {\bf (approximate tangential normality in reflexive spaces).} Let $\O$ be a
subset of a reflexive space $X$, and let $\ox\in\O$. Then given
any $d\in\Lm_w=T(\ox;\O)$ and $x^*\in\Hat N(d;\Lm_w)$, we have
\eqref{eq:(AN)} whenever sequences $\{x_k\}\in\mT^w_d$ and $t_k\dn
0$ are taken from the construction of $\mT^w_d$. In particular,
the set $\O$ enjoys the TAN property at $\ox$.
\end{Proposition}
{\bf Proof.} Assume that $\ox=0$ for simplicity. Pick any $\ve>0$ and by the definition of Fr\'echet normals find
$\dd\in(0,\ve)$ such that
\begin{equation}\label{eq:Lem(AN)}
\la x^*,v-d\ra\le\frac{\ve}{2}\|v-d\|\;\mbox{ for all }\;v\in\Lm_w\cap(d+\dd\B).
\end{equation}
Fix any sequences $\{x_k\}\in\mT^w_d$ and $t_k\dn 0$ from the
formulation of the proposition and show that property
\eqref{eq:(AN)} holds with the numbers $\ve$ and $\dd$ chosen
above. Supposing the contrary, find $\{x_k\}\in\mT^w_d$ and the
corresponding sequence $t_k\dn0$ such that
\begin{equation*}
\lim_{k\to\infty}\Big\{\sup\frac{\la x^*,z-x_k\ra}{t_k}\Big|\;z\in\O\cap(B(x_k+t_k\dd\B)\Big\}>2\ve\dd
\end{equation*}
along some subsequence of $k\in\N$, with no relabeling here and in
what follows. Hence there is a sequence of
$z_k\in\cap(x_k+t_k\dd\B)$ along which
\begin{equation*}
\frac{\la x^*,z_k-x_k\ra}{t_k}>\ve\dd\;\mbox{ for }\;k\in\N.
\end{equation*}
Taking into account the relationships
\begin{equation*}
\Big\|\frac{z_k}{t_k}-\frac{x_k}{t_k}\Big\|\le\dd\;\mbox{ and }\;
\frac{x_k}{t_k}\st{w}{\lto}d\;\mbox{ as }\;k\to\infty,
\end{equation*}
we get that the sequence $\disp\Big\{\frac{x_k}{t_k}\Big\}$ is bounded in $X$, and so
is $\Big\{\disp\frac{z_k}{t_k}\Big\}$. Since any bounded sequence in a
reflexive Banach space contains a weakly convergent subsequence, we may
assume with no loss of generality that the sequence $\disp\Big\{\frac{z_k}{t_k}\Big\}$ weakly converges to some
$v\in X$ as $k\to\infty$. It follows from the weak convergence of this sequence that
\begin{equation*}
\|v-d\|\le\liminf_{k\to\infty}\Big\|\frac{z_k}{t_k}-\frac{x_k}{t_k}\Big\|\le\dd.
\end{equation*}
This allows us to conclude that
\begin{equation*}
\la x^*,v-d\ra\ge\ve\dd>\frac{\ve}{2}\dd\ge\frac{\ve}{2}\|v-d\|,
\end{equation*}
which contradicts (\ref{eq:Lem(AN)}) and thus completes the proof of the proposition. $\h$\vspace*{0.05in}

The next theorem is the main result of this section showing that
the TAN property of a closed set in an Asplund space implies the
TNE property of the weak contingent cone to this set at the
reference point. This unconditionally justifies the latter
property in reflexive spaces.

\begin{Theorem}\label{Thm:TNE-Asplund} {\bf (TNE property in Asplund spaces).} Let $\O$ be a closed subset of an
Asplund space $X$, and let $\ox\in\O$. Assume that $\O$ has the
tangential approximate normality property at $\ox$. Then the weak
contingent cone $\Lm_w=T(\ox;\O)$ is tangentially normally
enclosed into $\O$ at this point. Furthermore, the latter TNE
property holds for any closed subset of a reflexive space.
\end{Theorem}
{\bf Proof.}  We are going show that the following holds in the Asplund space setting under the TAN property of $\O$ at
$\ox$:
\begin{equation}\label{an1}
\Hat N(d;\Lm_w)\subset N(\ox;\O)\ \mbox{ for all }\ d\in\Lm,\;\|d\|=1,
\end{equation}
which is obviously equivalent to $N(0;\Lm_w)\subset N(\ox;\O)$,
the TNE property of the weak contingent cone $\Lm_w$. Then the
second conclusion of the theorem in reflexive spaces immediately
follows from Proposition~\ref{Lem:(AN)}. Assume without loss of
generality that $\ox=0$.

To justify \eqref{an1}, fix $d\in\Lm_w$ and $x^*\in\Hat
N(d;\Lm_w)$ with $\|d\|=1$ and $\|x^*\|=1$. Taking
$\{x_k\}\in\mT^w_d$ from Definition~\ref{Def:(AN)}, it follows
that for any $\ve$ there is $\dd<\ve$ such that \eqref{eq:(AN)}
holds with $\ox=0$. Hence
\begin{equation}\label{an2}
\la x^*,z-x_k\ra\le 3t_k\ve\dd\;\mbox{ whenever }\;z\in Q:=\O\cap(x_k+t_k\dd\B),\quad k\in\N.
\end{equation}
Consider further the function
\begin{equation*}
\ph(z):=-\la x^*,z-x_k\ra,\quad z\in Q,
\end{equation*}
for which we have by \eqref{an2} that
\begin{equation*}
\ph(x_k)=0\le\inf_{z\in Q}\ph(z)+3t_k\ve\dd.
\end{equation*}
Setting $\lm:=\frac{t_k\dd}{3}$ and $\Tilde\ve:=3t_k\ve\dd$, we
apply the Ekeland variational principle (see, e.g.,
\cite[Theorem~2.26]{m-book1}) with $\lm$ and $\Tilde\ve$ to the
function $\ph$ on $Q$. In this way we find $\Tilde x\in Q$ such
that $\|\Tilde x-x_k\|\le\lm$ and $\Tilde x$ minimizes the
perturbed function
\begin{equation*}
\psi(z):=-\la x^*,z-x_k\ra+\frac{\Tilde\ve}{\lm}\|z-\Tilde
x\|=-\la x^*,z-x_k\ra+9\ve\|z-\Tilde x\|,\quad z\in Q.
\end{equation*}
Applying now the generalized Fermat rule to $\psi$ at $\Tilde x_k$
and then the fuzzy sum rule in the Asplund space setting (see,
e.g., \cite[Lemma~2.32]{m-book1}) gives us
\begin{equation}\label{an3}
0\in-x^*+(9\ve+\lm)\B^*+\Hat N(\Tilde x_k;Q)
\end{equation}
with some $\Tilde x_k\in\O\cap(\Tilde x+\lm\B)$. The latter means that
\begin{equation*}
\|\Tilde x_k-x_k\|\le\|\Tilde x_k-\Tilde x\|+\|\Tilde x-x_k\|\le2\lm<t_k\dd.
\end{equation*}
Hence $\Tilde x_k$ belongs to the interior of the ball centered at $\Tilde x$ with radius $t_k\dd$, which
implies that $\Hat N(\Tilde x_k;Q)=\Hat N(\Tilde x_k;\O)$. Thus we get from \eqref{an3} that
\begin{equation*}
x^*\in\Hat N(\Tilde x_k;\O)+(9\ve+\lm)\B^*,\quad k\in\N.
\end{equation*}
Letting there $k\to\infty$ and then $\ve\dn 0$ gives us $\Tilde
x_k\to\ox$ and $x^*\in N(\ox;\O)$. This justifies \eqref{an1} and
 completes the proof of the theorem. $\h$

\begin{Corollary}\label{Cor:TNE-Rn} {\bf (TNE property of the contingent cone in finite dimensions).} Let a set
$\O\subset\R^n$ be closed  around $\ox\in\O$. Then the contingent
cone $T(\ox;\O)$ to $\O$ at $\ox$ is tangentially normally
enclosed into $\O$ at this point, i.e., we have
\begin{equation}\label{an4}
N(0;\Lm)\subset N(\ox;\O)\;\mbox{ with }\;\Lm:=T(\ox;\O).
\end{equation}
\end{Corollary}
{\bf Proof}. It follows from Theorem~\ref{Thm:TNE-Asplund} due to $T(\ox;\O)=T_w(\ox;\O)$ in $\R^n$. $\h$\vspace*{0.05in}

Note that another proof of inclusion \eqref{an4} in $\R^n$ can be found in \cite[Theorem~6.27]{Rockafellar-Wets-VA}.

\section{Contingent and Weak Contingent Extremal Principles for Countable and Finite Systems of Closed Sets}
\setcounter{equation}{0}

By {\em tangential extremal principles} we understand results
justifying the validity of extremality conditions defined in
Section~3 for countable and/or finite systems of closed sets at
the corresponding {\em tangential extremal points}. Note that,
given a system of $\Lm=\{\Lm_i\}$-approximating cones to a set
system $\{\O_i\}$ at $\ox$, the results ensuring the fulfillment
of the $\Lm$-tangential extremality conditions at $\Lm$-tangential
local extremal points are directly induced by an appropriate conic
extremal principle applied to the cone system $\{\Lm_i\}$ at the
origin. It is remarkable, however, that for {\em tangentially
normally enclosed} cones $\{\Lm_i\}$ we simultaneously ensure the
fulfillment of the {\em limiting extremality conditions} of
Definition~\ref{Def:EPinf}(c) at the corresponding tangential
extremal points. As shown in Section~6, this is the case of the
contingent cone in finite dimensions and of the weak contingent
cone in reflexive (and also in Asplund) spaces.

In this section we pay the main attention to deriving the {\em
contingent} and {\em weak contingent extremal principle} involving
the aforementioned extremality conditions for countable and finite
systems of sets and finite-dimensional and infinite-dimensional
spaces. Observe that in the case of countable collections of sets
the results obtained are the first in the literature, while in the
case of finite systems of sets they are independent of the those
known before being applied to different notions of tangential
extremal points; see the discussions in Section~3.\vspace*{0.05in}

We begin with the contingent extremal principle for countable
systems of arbitrary closed sets in finite-dimensional spaces.

\begin{Theorem}\label{Thm:TEPinf} {\bf (contingent extremal principle for countable sets systems in finite dimensions).}
Let $\ox\in\bigcap_{i=1}^\infty\O_i$ be a contingent local
extremal point of a countable system of closed sets
$\{\O_i\}_{i\in\N}$ in $\R^n$. Assume that the contingent cones
$T(\ox;\O_i)$ to $\O_i$ at $\ox$ are nonoverlapping
\begin{equation*}
\bigcap_{i=1}^\infty\Big\{T(\ox;\O_i)\Big\}=\big\{0\big\}.
\end{equation*}
Then there are normal vectors
\begin{equation*}
x^*_i\in N(0;\Lm_i)\subset N(\ox;\O_i)\;\mbox{ for }\;\Lm_i:=T(\ox;\O_i)\;\mbox{ as }\;i\in\N
\end{equation*}
satisfying the extremality conditions in \eqref{eq:EPinf}.
\end{Theorem}
{\bf Proof}. This result follows from combining Theorem~\ref{Thm:CEP} and
Corollary~\ref{Cor:TNE-Rn}. $\h$\vspace*{0.05in}

Consider further systems of finitely many sets
$\{\O_1,\ldots,\O_m\}$ in Asplund spaces and derive for them the
weak contingent extremal principle. Recall that a set $\O\subset
X$ is {\em sequentially normally compact} (SNC) at $\ox\in\O$ if
for any sequence $\{(x_k,x^*_k)\}_{k\in\N}\subset\O\times X^*$ we
have the implication
\begin{equation*}
\big[x_k\to\ox,\;x^*_k\st{w^*}{\to}0\;\mbox{ with }\;x^*_k\in\Hat N(x_k;\O),\;k\in\N\big]\Longrightarrow\|x^*_k\|\to 0\;
\mbox{ as }\;k\to\infty.
\end{equation*}
In \cite[Subsection~1.1.4]{m-book1}, the reader can find a number
of efficient conditions ensuring the SNC property, which holds in
rather broad infinite-dimensional settings. The next proposition
shows that the SNC property of TAN sets is inherent by their weak
contingent cones.

\begin{Proposition}\label{Prop:SNCwB} {\bf (SNC property of weak contingent cones).} Let $\O$ be a closed subset of an
Asplund space $X$ satisfying the tangential approximate normality
property at $\ox\in\O$. Then the weak contingent cone
$T_w(\ox;\O)$ is SNC at the origin provided that $\O$ is SNC at
$\ox$. In particular, in reflexive spaces the SNC property of a
closed subset $\O$ at $\ox$ unconditionally implies the SNC
property of its weak contingent cone $T_w(\ox;\O)$ at the origin.
\end{Proposition}
{\bf Proof.} To justify the SNC property of $\Lm_w:=T_w(\ox;\O)$
at the origin, take sequences $d_k\to 0$ and $x^*_k \in\Hat
N(d_k;\Lm_w)$ satisfying $x^*_k\st{w^*}{\to}0$ as $k\to\infty$.
Using the TAN property of $\O$ at $\ox$ and following the proof of
Theorem~\ref{Thm:TNE-Asplund}, we find sequences $\ve_k\dn0$ and
$\Tilde x_k\st{\O}{\to}\ox$ such that
\begin{equation*}
x^*_k\in\Hat N(\Tilde x_k;\O)+\ve_k\B^*\;\mbox{ for all }\;k\in\N.
\end{equation*}
Hence there are $\Tilde x^*_k\in\Hat N(\Tilde x_k;\O)$ with
$\|\Tilde x^*_k-x^*_k\|\le\ve_k$, which implies that $\Tilde
x^*_k\st{w^*}{\to}0$ as $k\to\infty$. By the SNC property of $\O$
at $\ox$ we get that $\|\Tilde x^*_k\|\to 0$, which yields in turn
that $\|x^*_k\|\to 0$ as $k\in\infty$. This justifies the SNC
property of $\Lm_w$ at the origin. The second assertion of this
proposition immediately follows from Proposition~\ref{Lem:(AN)}.
$\h$\vspace*{0.05in}

Now we are ready to establish the weak contingent extremal
principle for systems of finitely many closed subsets of Asplund
spaces in both approximate and exact forms.

\begin{Theorem}\label{Thm:TEP-Asp}{\bf (weak contingent extremal principle for finite systems of sets in Asplund spaces).} Let
$\ox\in\bigcap_{i=1}^m\O_i$ be a weak contingent local extremal
point of the system $\{\O_1,\ldots,\O_m\}$ of closed sets in an
Asplund space $X$. Assume that all the sets $\O_i$,
$i=1,\ldots,m$, have the TAN property at $\ox$, which is automatic
in reflexive spaces. Then the following versions of the weak
contingent extremal principle hold:

{\bf (i)} {\sc Approximate version}: for any $\ve>0$ there are $x^*_i\in N(\ox;\O_i)$ as $i=1,\ldots,m$ satisfying
\begin{equation}\label{eq:TEP-A}
\|x^*_1+\ldots+x^*_m\|\le\ve\;\mbox{ and }\;\quad\|x^*_1\|+\ldots+\|x^*_m\|=1.
\end{equation}

{\bf (ii)} {\sc Exact version}: if in addition all but one of the
sets $\O_i$ as $i=1,\ldots,m$ are SNC at $\ox$, then there exist
$x^*_i\in N(\ox;\O_i)$ as $i=1,\ldots,m$ satisfying
\begin{equation}\label{asp2}
x^*_1+\ldots+x^*_m=0\;\mbox{ and }\;\|x^*_1\|+\ldots+\|x^*_m\|=1.
\end{equation}
\end{Theorem}
{\bf Proof}. It follows from Proposition~\ref{cone-eq} that the
cone system $\{\Lm^i_w=T_w(\ox;\O_i)\}$ as $i=1,\ldots,m$ is
extremal at the origin in the conventional sense \eqref{eq:ESfin}.
Applying to it the approximate extremal principle from
\cite[Theorem~2.20]{m-book1}, for any $\ve>0$ we find
$x_i\in\Lm_w^i$ and $x^*_i\in\Hat N(x_i;\Lm_w^i)$ as $i=1,...,m$
such that all the relationships in (\ref{eq:TEP-A}) hold. Then
\begin{equation*}
x^*_i\in\Hat N(x_i;\Lm^i_w)\subset N(0;\Lm^i_w)\subset N(\ox;\O_i),\quad i=1,\ldots,m,
\end{equation*}
by Proposition~\ref{Prop:N(Lm)} and Theorem~\ref{Thm:TNE-Asplund}, which justifies assertion (i).

Now to justify (ii), observe that all but one of the cones
$\Lm^i_w$ are SNC at the origin by Proposition~\ref{Prop:SNCwB}.
Thus the conclusion of (ii) follows from the exact extremal
principle in \cite[Theorem~2.22]{m-book1} and
Theorem~\ref{Thm:TNE-Asplund} established above. $\h$

\end{document}